\documentclass[a4paper, 10pt, conference]{ieeeconf}      

\usepackage[dvips]{epsfig}
\usepackage{amssymb,graphics,amsmath,amsthm,amsopn,amstext,amsfonts}

\newtheorem{proposition}{Proposition}
\newtheorem{theorem}{Theorem}

\newtheorem{lemma}{Lemma}

\newcommand{\be}{\begin{equation}}
\newcommand{\ee}{\end{equation}}
\newcommand{\bea}{\begin{eqnarray}}
\newcommand{\eea}{\end{eqnarray}}
\newcommand{\beas}{\begin{eqnarray*}}
\newcommand{\eeas}{\end{eqnarray*}}
\newcommand{\nn}{\nonumber}
\newcommand{\bbm}{\begin{bmatrix}}
\newcommand{\ebm}{\end{bmatrix}}
\newcommand{\matl}{\left[ \begin{array}}
\newcommand{\matr}{\end{array} \right]}
\newcommand{\la}{\label}
\newcommand{\La}{\mathcal{L}}

\newcommand{\cJ}{\mathcal{J}}
\newcommand{\bbr}{\mathbb{R}}
\newcommand{\mfrak}{\mathfrak}
\newcommand{\Z}{\mathbb{Z}}
\newcommand{\sumk}{\sum_{k=1}^{N}}
\newcommand{\Ta}{\mathtt{T}}
\newcommand{\T}{^{\mbox{\small T}}}
\newcommand{\mbf}{\mathbf}
\newcommand{\tr}{\mbox{trace}}
\newcommand{\di}{\mathrm{d}}
\newcommand{\wh}{\widehat}
\newcommand{\wt}{\widetilde}
\newcommand{\lan}{\langle}
\newcommand{\ran}{\rangle}
\newcommand{\Lam}{\Lambda}

\begin{document}




\title {Optimal Attitude Estimation and Filtering Without Using Local Coordinates \\
Part I: Uncontrolled and Deterministic Attitude Dynamics} 

\author{Amit K. Sanyal \\ 
Department of Mechanical and Aerospace Engineering, \\
Arizona State University, \\
Tempe, AZ 85287 \\
\tt{ sanyal@asu.edu }} 

\maketitle
\pagestyle{empty}
\renewcommand{\baselinestretch}{0.8}
\begin{abstract}
There are several attitude estimation algorithms in existence, all of 
which use local coordinate representations for the group of rigid body 
orientations. All local coordinate representations of the group of orientations 
have associated problems. While minimal coordinate representations exhibit 
kinematic singularities for large rotations, the quaternion representation 
requires satisfaction of an extra constraint. This paper treats the attitude 
estimation and filtering problem as an optimization problem, without using 
any local coordinates for the group of rotations. An attitude determination 
algorithm and attitude estimation filters are developed, that minimize the 
attitude and angular velocity estimation errors. For filter propagation, the 
attitude kinematics and deterministic dynamics equations (Euler's equations) 
for a rigid body in an attitude dependent potential are used. Vector attitude 
measurements are used for attitude and angular velocity estimation, with or 
without angular velocity measurements.  
\end{abstract}

\section{Introduction}

Attitude estimation is often a prerequisite for controlling aerospace and 
underwater vehicles, mobile robots, and other mechanical systems moving in 
space. Hence, attitude estimation of a rigid body has applications in spacecraft 
and aircraft dynamics, unmanned vehicle dynamics, and robot dynamics, including 
walking robots. While attitude sensors and the control tasks for which attitude 
feedback are required may be different in these different applications, the 
fundamental importance of obtaining accurate attitude data remains common to 
all these applications. In this paper, a new look at the attitude estimation 
problem is provided, which has two essentially new features: (1) the attitude 
is globally represented without using any local coordinates and the nonlinear 
attitude dynamics equation (Euler's equation) for rigid bodies is used, and (2) 
the filter obtained is not a Kalman or extended Kalman filter. A global attitude 
representation has been recently used for partial attitude estimation with a 
linear dynamics model (see \cite{rehu}). However, to the author's knowledge, 
total attitude estimation using a global attitude representation and a full 
nonlinear attitude kinematics and dynamics model (without linearization) has not 
been done before. 

Spacecraft attitude determination and filtering is perhaps the oldest application
for attitude estimation algorithms, and the attitude determination problem for a 
spacecraft from vector measurements was first posed in \cite{wah}. A sample of 
the literature in this area can be found in \cite{bosh, cram, mark, shus1, shus2, wah}. 
Applications of attitude estimation to unmanned vehicles and robots can be found in 
\cite{badu, rehu, rosb, vafo}. Algorithms that are typically used for attitude 
estimation in such applications are based upon local coordinate representations 
of the group of rotations, like quaternions, Rodrigues parameters, or Euler 
angles. As is well known, minimal coordinate representations of the rotation group, 
like Euler angles, Rodrigues parameters, and modified Rodrigues parameters (see 
\cite{cram2}), usually lead to geometric or kinematic singularities. Quaternion 
representation of the attitude matrix is commonly used, particularly in spacecraft 
applications, where the quaternion estimation (QUEST) algorithm and its several 
variants have been in use for quite some time (\cite{bosh, shus2, psia}). Besides  
the extra constraint (of unit norm) that one needs to impose on the quaternion, the 
quaternion representation for a given rotation depends on the sense of rotation 
used to define the principal angle, and hence can be defined in one of two ways. 
Local coordinate representations of the attitude usually lead to use of the extended 
Kalman filter (EKF) as an estimator for attitude and angular velocity. It is well 
known that the EKF has problems with convergence and stability in the case of large 
initial condition errors \cite{crajun}. The attitude determination algorithm presented 
here does not use any local coordinate representation of the attitude, and is hence 
free of the drawbacks associated with such local representations. Nonlinear attitude 
estimation filters for a rigid body in an attitude-dependent potential field are 
also developed using this attitude determination algorithm. These are optimal 
nonlinear filters that minimize the attitude and angular velocity estimation errors 
at each measurement instant, and are hence free of the stability issues confronting 
extended Kalman filters.  

A brief outline of this paper is given here. In Section \ref{sec:atder}, the 
attitude determination problem for vector measurements with measurement noise is 
introduced, and a global attitude determination algorithm which gives a global 
minimum of the attitude estimation error is presented. This section also presents 
some simulation results that demonstrate the applicability of this attitude 
determination algorithm. Section \ref{sec:atfil} introduces an attitude dynamics 
model for a free rigid body in a potential field, where the inertia properties of 
the body are assumed to be perfectly known. This deterministic dynamics model is 
used to create a filter that estimates both the attitude and the angular velocity. 
Two cases are considered here: the presence of angular velocity measurements, 
and their absence; and the filter algorithms for both cases are presented. 
The paper is concluded in Section \ref{sec:conc} with a summary of results presented, 
and a discussion on future enhancements to the filter algorithm developed here. 

\section{Attitude Determination Problem}\la{sec:atder}

The attitude of a rigid body in space is a representation of the orientation of 
a body-fixed coordinate frame to an inertial frame. The principal axes of the 
body and inertial frames are related by a linear transformation given by a 
proper orthogonal matrix, which is usually referred to as the rotation or 
orientation matrix. The rotation matrix may be represented by various sets of 
coordinates, like the Euler angles, quaternions, or Rodrigues parameters (see 
\cite{gold, green}). The rotation matrices are proper (determinant=+1) orthogonal 
matrices that form a group under matrix multiplication; this abstract group 
is denoted $\mbf{SO}(3)$. Hence, the group $\mbf{SO}(3)$ is the compact Lie 
group of orientation-preserving isometries on $\bbr^3$, and we represent it 
using the set of $3\times 3$ proper orthogonal matrices,  
\[ C\in\bbr^{3\times 3},\ \mbox{ s. t. }\ C\T C=I_3= CC\T,\;\ \det (C)=1. \]

\subsection{Attitude Determination from Vector Measurements}

We now formulate the attitude determination problem from vector measurements. 
Let the direction vectors of a few known points in an inertial frame 
$\mathcal{I}$ for $\bbr^3$ be given by 
\[ e_i,\;\ i=1,2,\ldots,n, \] 
and their corresponding direction vectors in a body-fixed frame 
$\mathcal{B}$ (fixed to a rigid body of interest) for $\bbr^3$ be
\[ b_i,\;\ i=1,\ldots,n. \] 
The inertial and body-fixed direction vectors are related by the rotation 
matrix $C$ which rotates the body frame into the inertial frame, such that 
\be e_i= C b_i\; \forall\ i\in \{1,2,\ldots,n\}. \la{eCb} \ee 
Note that the convention followed for the rotation matrix $C$ in (\ref{eCb}) 
is that in \cite{bbcm, mara}, while the reverse convention of a rotation
matrix taking the inertial frame to the body frame is used in most of the 
other literature cited here. The convention used here makes it easier to 
represent the body kinematics and dynamics in the body frame and the equations 
of motion are {\em left-invariant}, i.e., invariant to left multiplication of 
$C$ by a non-singular matrix. The direction vectors $b_i$ when measured from the 
body (e.g., a spacecraft), usually contain additive measurement errors and 
the measured direction vectors may not coincide with the actual $b_i=C\T e_i$. 
Let the measured direction vectors be given by
\[ \wt{b}_i=b_i+\nu_i, \]
where $\nu_i$ are measurement errors that are usually assumed to be 
Gaussian with zero mean. 

The attitude determination problem consists of finding an estimate $\wh{C}$ 
of the rotation matrix $C$ such that the errors 
\[ e_i- \wh{C}\wt{b}_i \]
are minimized. The least squares attitude estimation problem would 
be to
\[ \mbox{Minimize }\frac12 \sum_{i=1}^n w_i(e_i-\wh{C}\wt{b}_i)\T (e_i-\wh{C}
\wt{b}_i) \]
with respect to $\wh{C}$ subject to $\wh{C}\in\mbf{SO}(3)$, where $w_i$ is 
a known wieght factor (positive) usually taken to correspond to the statistical 
standard deviation of the $i$th measured vector. This problem is also 
known as Wahba's problem \cite{wah}. In this work, the weight factors are 
considered as design parameters. We define the $3\times n$ matrices 
\[ E=[e_1\;\ e_2\;\ \ldots\;\ e_n],\;\,\ \wt{B}=[\wt{b}_1\;\ \wt{b}_2\;\ 
\ldots\;\ \wt{b}_n]. \] 
The assumption here and throughout the rest of this paper is that 
both $E$ and $\wt{B}$ are of rank 3; otherwise the attitude determination 
problem is ill-posed. We introduce the trace inner product on the space 
of real $n_1\times n_2$ matrices, i.e., if $A_1,\ A_2\in\bbr^{n_1\times n_2}$, 
then 
\[ \langle A_1,A_2\rangle= \tr (A_1\T A_2). \]
The above attitude determination problem can then be restated as follows:
\be \mbox{Minimize }\cJ_0= \frac12 \langle E-\wh{C}\wt{B}, (E-\wh{C}\wt{B})
W\rangle,\;\ \wh{C}\in \mbf{SO}(3), \la{atprob} \ee 
where the estimate of the rotation matrix (attitude) $\wh{C}$ is the only 
unknown and $W=\mbox{diag}(w_i)$ is the positive diagonal weight matrix. We 
can extemize the cost function $\cJ_0$ by taking the first variation with 
respect to $\wh{C}$ and setting it to zero since $\wh{C}$ is the only 
unknown to be determined in this problem. The extremal solution to this 
problem is given by 
\bea \delta_{\wh{C}}\cJ_0&=&\frac12\langle -\delta\wh{C}\wt{B},(E-\wh{C}
\wt{B})W\rangle+\frac12\langle E-\wh{C}\wt{B}, \nn \\ -\delta\wh{C}\wt{B}W
\rangle &=& \langle (E-\wh{C}\wt{B})W\wt{B}\T,-\delta\wh{C}\rangle \nn \\
&=& -\tr(\wt{B}W(E-\wh{C}\wt{B})\T\delta\wh{C})=0, \la{solp} \eea
where $\delta\wh{C}$ is a variation in $\wh{C}\in\mbf{SO}(3)$. Since $\delta
\wh{C}$ is in $T_{\wh{C}}\mbf{SO}(3)$, the tangent space to $\mbf{SO}(3)$ 
at $\wh{C}$, it has the form 
\be \delta\wh{C}= \wh{C}U,\;\ U\in \mfrak{so}(3). \la{delC} \ee
Hence, from (\ref{solp}) and (\ref{delC}), we get 
\bea \tr(\wt{B}W(E-\wh{C}\wt{B})\T\wh{C}U)=0 \Leftrightarrow \wt{B}W(E\T 
\wh{C}-\wt{B}\T) \nn \\ \mbox{is symmetric} \Leftrightarrow \wt{B}WE\T\wh{C}\ 
\mbox{is symmetric}, \la{BETC} \eea
since $U$ is skew-symmetric. The above result can be recast into 
the following form:
\be L\T\wh{C}= \wh{C}\T L,\;\ L=EW\wt{B}\T, \la{LTC} \ee
and $L$ is known since $E$ is known, $\wt{B}$ is known from measurements, and 
$W$ is known as a design parameter. 

The following result gives a necessary condition for the attitude matrix $C$ 
that satisfies (\ref{solp}), and is equivalent to equation (\ref{LTC}). 
\begin{lemma}
Define the linear map $M_L : \mbf{SO}(3)\rightarrow\mfrak{so}(3)$ by 
\be M_L(C)=C\T L-L\T C,\;\ C\in\mbf{SO}(3), \la{MLdef} \ee
where $L$ is as defined by (\ref{LTC}). If $C\in\mbf{SO}(3)$ is in the 
kernel of this map, then $C$ is of the form
\be C=SL,\;\ S=S\T, \la{CSL} \ee
i.e., $C=SL$ where $S$ is a $3\times 3$ symmetric matrix.
\la{CSLlem}
\end{lemma}
\noindent {\em Proof}: If $C$ is in the kernel of $M_L$, then 
\[ C\T L=L\T C\ \Rightarrow\ LC\T= CL\T. \]
Hence, $D=LC\T$ is symmetric and we could express $L=DC$, where $D=D\T$ 
is symmetric. Now from our earlier assumptions, $L=EWB\T$ is non-singular, 
since $E$ and $B$ are of rank 3, and $W\in\bbr^{n\times n}$ is positive 
definite. Thus $D=LC\T$ is also non-singular. Hence, we get 
\[ DC=L \ \Rightarrow\ C=D^{-1}L=SL, \]
where $S=D^{-1}$ is symmetric. This proves the result. \qed \\
This result is a special case of Proposition 1 in \cite{san}, in which 
$C$ is replaced by a matrix whose row vectors form an orthonormal set. 
However, the above result does not give the unique solution to the 
attitude determination problem (\ref{atprob}) since it does not give 
an expression for $S$, from which the estimate $\wh{C}$ of the unknown 
attitude can be determined. 

To obtain the $\wh{C}$ that minimizes the 
cost function $\cJ_0$, we apply the sufficient condition for a minimum 
by taking its second variation with respect to $\wh{C}$. The first 
variation of $\cJ_0$ in (\ref{solp}) can be written as 
\[ \delta_{\wh{C}}\cJ_0 = \lan L,-\delta\wh{C}\ran= -\lan L,\wh{C}U\ran. \]
Thus, a sufficient condition for $\wh{C}$ to minimize the cost function 
$\cJ_0$ is as follows:
\bea \delta^2_{\wh{C}}\cJ_0 &=& -\tr\big(L\T\wh{C}U^2\big)-\tr\big(L\T\wh{C}
\delta U\big) \nn \\ &=& -\tr\big(L\T\wh{C}U^2\big) >0. \la{secn} \eea 
This condition, along with Lemma \ref{CSLlem}, leads to the following result. 
\begin{proposition}
The cost function $\cJ_0$ in (\ref{atprob}) is minimized by $\wh{C}=SL$ such 
that the symmetric matrix $S$ is positive definite. 
\la{promin} 
\end{proposition}
\noindent {\em Proof}: From Lemma \ref{CSLlem}, we know that a necessary 
condition for the minimizing $\wh{C}$ would be $\wh{C}=SL$, where $S$ is a 
symmetric matrix. Hence 
\[ D= L\T\wh{C}= L\T SL \]
is symmetric. From condition (\ref{secn}), we have $\tr\big(DU^2\big)<0$. 
Since $U\in\mfrak{so}(3)$, $U^2$ is symmetric and has negative definite 
trace. Let $Q_1, Q_2\in\mbf{O}(3)$ be such that 
\[ D=Q_1\Lam_1 Q_1\T,\;\ -U^2= Q_2\Lam_2 Q_2\T, \]
are the spectral decompositions of $D$ and $-U^2$ respectively, and $\Lam_1$, 
$\Lam_2$ are diagonal. Then 
\beas c= \tr(-DU^2)= \tr\big(Q_1\Lam_1 Q_1\T Q_2\Lam_2 Q_2\T\big) &>& 0 \\
\Rightarrow \tr\big(\Lam_1 Q_1\T Q_2\Lam_2 Q_2\T Q_1\big) &>& 0 \\ 
\Rightarrow \tr\big(\Lam_1 P\Lam_2 P\T\big) &>& 0 \eeas where 
$P=Q_1\T Q_2\in \mbf{O}(3)$. We denote the columns of $P$ by $p_1,\ p_2,\ p_3
\in\bbr^3$ and the entries of the diagonal matrix $\Lam_1$ as $l_1$, $l_2$, 
and $l_3$. Since $U\in\mfrak{so}(3)$, the eigenvalues of $-U^2$ (the 
entries of the diagonal matrix $\Lam_2$) are 0, $u_1>0$, and $u_2>0$. 
Thus, the sufficient condition (\ref{secn}) is equivalent to 
\beas c= u_1\tr(\Lam_1 p_2 p_2\T)+ u_2\tr(\Lam_1 p_3 p_3\T) &>& 0 \\
\Rightarrow (l_1 p_{21}^2+ l_2 p_{22}^2+ l_3 p_{23}^2) &>& 0\\ \mbox{ and }
\ (l_1 p_{31}^2+ l_2 p_{32}^2+ l_3p_{33}^2) &>& 0, \eeas 
which is possible for arbitrary $p_2,\ p_3\in\bbr^3$ if and only if 
$l_1>0$, $l_2>0$, and $l_3>0$. Hence, the sufficient conditon (\ref{secn}) 
is equivalent to $D=Q_1\Lam_1 Q_1\T$ being positive definite, i.e., 
$v\T Dv >0$ for any $v\in\bbr^3$. This in turn implies that $S=(L\T)^{-1}
DL^{-1}$ is also positive definite, since $u\T Su= v\T Dv>0$ where 
$u=Lv$. \qed 

Lemma \ref{CSLlem} and Proposition \ref{promin} give necessary and sufficient 
conditions, respectively, for the attitude matrix that minimizes $\cJ_0$. 
The following result gives an unique attitude matrix $\wh{C}\in\mbf{SO}(3)$ 
that solves the attitude determination problem (\ref{atprob}) and satisfies 
equations (\ref{LTC}) and (\ref{secn}).  
\begin{theorem}
The unique minimizing solution to the attitude determination problem 
(\ref{atprob}) is given by 
\be \wh{C}= SL,\;\ S=Q\sqrt{(RR\T)^{-1}}Q\T, \la{soln} \ee
where 
\be L=QR,\;\ Q\in\mbf{SO}(3), \la{LQR} \ee
and $R$ is upper triangular and non-singular (since $L=EWB\T$ is 
non-singular); this is the QR decomposition of $L$. The matrix square root 
used here is the positive definite (principal) square root of a positive 
definite symmetric matrix. 
\la{atde}
\end{theorem}
\noindent {\em Proof:} From Lemma \ref{CSLlem}, we know that $\wh{C}=SL=
SEW\wt{B}\T$ where $S=S\T$, is a necessary condition for the extremal 
solution. From Proposition \ref{promin}, the equivalent condition to the 
sufficient condition (\ref{secn}) is that the symmetric matrix $S$ has to 
be positive definite. Using the QR decomposition of $L$ given by (\ref{LQR}), 
we can express the orthogonality condition of $\wh{C}$ as follows:  
\[ \wh{C}\wh{C}\T= SQRR\T Q\T S= I_3. \]
Since $S$ is symmetric, $S$ is given by 
\[ S=\sqrt{Q(RR\T)^{-1}Q\T}= Q\sqrt{(RR\T)^{-1}}Q\T, \]
where the principal (positive definite) square root is taken, as given by 
equation (\ref{soln}). This makes $S$ positive definite as well. By 
construction, $\wh{C}=SL$ satisfies $\wh{C}\wh{C}\T= I_3$. 
Now we check the determinant of $\wh{C}=SL$, as follows: 
\beas \det \wh{C}&=& \det S \det L \\
&=& (\det Q\T \det \sqrt{(RR\T)^{-1}}\det Q) (\det Q\det R) \\
&=& \frac{(\det Q)^2}{\sqrt{\det R \det R\T}} \det Q\det R \\
&=& \frac1{\sqrt{(\det R)^2}} \det Q \det R \\
&=& \frac1{\det R} \det Q \det R= \det Q= 1, \eeas
since $Q\in\mbf{SO}(3)$. This proves that $\wh{C}\in\mbf{SO}(3)$, and is hence 
the unique minimal solution to the attitude determination problem (\ref{atprob}). 
\qed \\
Although we have used the QR decomposition for the matrix $L$ here, one can 
use the singular value decomposition or any other decomposition using orthogonal 
matrices, to show this result. 

We next show that the attitude estimate given 
by this algorithm is unbiased, i.e., in the absence of measurement errors,  
this estimate gives the actual attitude. 
\begin{proposition}
The attitude determination algorithm given by equations (\ref{soln})-(\ref{LQR}) 
gives an unbiased estimate of the attitude. 
\la{unbias} 
\end{proposition} 
\noindent {\em Proof:} Let us assume that there is no error in the measurement 
of body vectors, i.e., $\wt{B}=B=[b_1\;\ b_2\;\ \ldots\;\ b_n]$. In that case, 
we have $E=CB$, and
\beas L &=& L_0 = EWB\T = EWE\T C \\
&\Rightarrow & C= (EWE\T)^{-1} L_0, \eeas
which, by Theorem \ref{atde}, implies that 
\[ S_0= Q_0\sqrt{(R_0 R_0\T)^{-1}}Q_0\T= (EWE\T)^{-1}, \]
where $L_0=Q_0 R_0$. This is equivalent to 
\beas EWE\T &=& Q_0\sqrt{(R_0 R_0\T)}Q_0\T, \\
\Leftrightarrow (EWE\T)^2 &=& Q_0 R_0 R_0\T Q_0\T. \eeas
But the right-hand side above is $L_0 L_0\T= (EWE\T C)(C\T EWE\T)$, which is 
equal to the left-hand side. Thus, we have $C=S_0 L_0$ where $S_0=(EWE\T)^{-1}$ 
as required. This proves that this algorithm is unbiased. \qed \\  
These results are used as the basis for attitude estimation filters obtained 
in Section \ref{sec:atfil}.

\subsection{Simulation results for Attitude Determination Algorithm}

We end this section with a simulated example of an attitude determination 
problem, where an attitude matrix is obtained from `measurements' of 
seven (unit) vectors, representing seven different directions in Euclidean 
3-space. The ``measured" vectors are given as normalized (unit) vectors. 
The simulation is carried out using a MATLAB program, which implements 
the attitude determination algorithm given in Theorem \ref{atde}. The simulated 
vectors ``measured" in the body frame have added Gaussian noise with a standard 
deviation of $\sigma=0.002\ \mbox{rads}\approx 0.115^\circ$ (which is 
relatively large compared to the capabilities of most modern attitude 
sensors). 

The data used in this simulation, in terms of the inertial unit vectors, 
and their ``measured" counterparts in the body frame, are as follows: \\
Simulation Data: 
\beas E=\matl{cccc} 0.3817 &0.3077& 0.2324& 0.3374 \\ -0.5450& 
-0.6045& -0.5824& -0.5675 \\ 0.7465& 0.7347& 0.7789& 0.7511 \end{array} \right. \\ 
\left. \begin{array}{ccc} 0.3161& 0.2975& 0.2807\\ -0.6582& -0.6046& -0.5912\\ 
0.6832& 0.7389& 0.7561 \matr, \\
\wt{B}=\matl{cccc} 0.1287& 0.0975& 0.1580& 0.1264 \\ -0.9628& 
-0.9843& -0.9833& -0.9750 \\ -0.2394& -0.1517& -0.0862& -0.1904 \end{array} 
\right. \\ 
\left. \begin{array}{ccc} 0.0210& 0.1020& 0.1249\\ -0.9904& -0.9829& -0.9836 
\\ -0.1414& -0.1404& -0.1279 \matr. \eeas 
Note that the vectors in the inertial frame are clustered together, and 
hence, so are the vectors in the body frame. This simulates direction 
vectors as would be measured by an optical instrument with a finite 
field of view, e.g., a star tracker. The ``actual" attitude matrix which 
takes the ``actual" body directions to the inertial directions, is 
assumed to be known for this simulation, and is given by
\be C=\matl{ccc} -0.2029& -0.1865& -0.9613\\ 0.6385& 0.7191& 
-0.2743\\ 0.7424& -0.6694& -0.0269 \matr. \ee 
As mentioned before, these simulated measurements of body directions 
correspond to added Gaussian noise of 0.002 radian to the ``actual" body 
directions. 

The results of this simulation, in the form of the attitude matrix 
determined by this algorithm $\wh{C}$, the error between the known ``actual" 
attitude matrix $C$ and the attitude $\wh{C}$, and the error $e=E-\wh{C}\wt{B}$, 
are given below. \\
Simulation Results: 
\beas \wh{C}= \matl{ccc} -0.2042& -0.1856& -0.9612\\ 
0.6386& 0.7190& -0.2745\\ 0.7420& -0.6698& -0.0283 \matr, \\
e_C=\wh{C}\T C-I_3=\matl{ccc} -0.0000& 0.0006& 0.0012\\ 
-0.0006& -0.0000& -0.0008\\ -0.0012& 0.0008& -0.0000 \matr, \\
E-\wh{C}\wt{B}=\matl{ccccccc} -0.0009& -0.0008& -0.0006& -0.0007 \\
-0.0007& -0.0008& -0.0000& 0.0005 \\ -0.0007& -0.0012& 0.0006& -0.0012 
\end{array} \right. \\ 
\left. \begin{array}{ccc} 0.0007& 0.0009& 0.0008\\ 
0.0003& 0.0009& 0.0009 \\ 0.0016& -0.0016& 0.0011 \matr. \eeas 
Note that the error in the attitude matrix is here specified as $\wh{C}\T C$ 
minus the identity matrix; one can also specify this error as $C-\wh{C}$. As 
defined here, the matrix $I_3+e_C=\wh{C}\T C\in\mbf{SO}(3)$, is the measure 
(in the group of rigid body rotations) of the attitude error. 
The maximum errors in these results are of the order of the measurement 
errors in the body vectors, which demonstrates the applicability of this 
attitude determination algorithm. 

\section{Attitude Estimation Filters for a Free Rigid Body in a Potential Field}
\la{sec:atfil} 

In this section, we develop an attitude estimation filter based on the 
attitude determination algorithm of Theorem \ref{atde} developed in 
the last section. We assume that the attitude dynamics is perfectly known, 
and is that of a free rigid body in a potential field, i.e., there are no 
applied (control) forces on the body. We leave the potential field to be 
general (could be uniform or central gravity, for example). Two cases are 
dealt with here: (1) the case without, and (2) the case with angular 
velocity measurements. Since we use the actual (continuous) nonlinear 
dynamics equations for filter propagation, the filters developed here 
are not extended Kalman filters; they are nonlinear filters. We also assume 
that the vector attitude measurements (and angular velocity measurements, 
if any), are made at discrete time instants. Hence, the filter equations 
obtained are of the continuous-discrete type. 

Let $\{t_k\}$, $k\in\Z$ non-negative, denote an increasing sequence of 
non-negative real numbers that coincide with time instants at which measurements 
of vectors in the body frame (and angular velocity measurments, if any) are 
taken. Let $\wt{B}_k\in\bbr^{3\times n_k}$ denote the set of $n_k$ vector 
measurments taken at time $t_k$ in the body frame; the columns of $\wt{B}_k$ 
denote the measured body vectors. If $E_k\in\bbr^{3\times n_k}$ denotes the 
vectors in the inertial frame, and $C_k$ is the actual attitude matrix from 
the body to the inertial frame, then 
\be \wt{B}_k= C_k\T E_k+ N_k, \la{bomeas} \ee
where the columns of $N_k\in\bbr^{3\times n_k}$ are the measurement errors 
in the body vectors. The measured body vectors are usually expressed as unit 
vectors. The angular velocity measurement at time  
$t_k$ is denoted by $\wt{\Omega}_k\in\mfrak{so}(3)$. At time $t_k$, $k\ge 1$, 
the attitude and angular velocity estimates obtained by propagating using 
the attitude kinematics and dynamics equations from time $t_{k-1}$ are denoted 
$\wh{C}_k^{-}$ and $\wh{\Omega}_k^-$ respectively, and the updated attitude 
at this time instant obtained from the attitude determination part of the 
filter (which is based on the algorithm of Section \ref{sec:atder}) is 
denoted $\wh{C_k}^+$. The angular velocity is also updated at the measurment 
time $t_k$, and the updated angular velocity estimate is denoted 
$\wh{\Omega}_k^+$. 

\subsection{Dynamics of Free Rigid Body in a Potential Field}

We first obtain the dynamics (equations of motion) of a free rigid body 
in a potential field in a compact geometric form, that is free of any 
particular coordinate description. 
The attitude kinematics is given by 
\be \dot C= C\Omega, \la{atkin} \ee
where $\Omega\in\mfrak{so}(3)$ denotes the angular velocity in the body 
frame. Let $\Lambda$ denote the symmetric positive definite inertia 
matrix of the rigid body. The Lagrangian for the rigid body in a potential 
field is given by 
\be \La (C,\Omega)=\frac12\langle \Omega,\Omega\Lambda\rangle- V(C), \la{Lagrb} \ee
where the first term is the kinetic energy, and $V(C)$ denotes the 
potential energy that is dependent on the attitude of the body. 

The equations of motion are obtained by applying Hamilton's principle to 
the action quantity 
\[ S= \int_0^T \La (C,\Omega) \di t, \]
and taking {\em reduced variations} on the group $\mbf{SO}(3)$ (see 
\cite{bbcm, mara}). The reduced variations at the point $(C,\Omega)\in\Ta 
\mbf{SO}(3)$ are given by (\cite{bbcm, mara}) 
\be \delta C= C\Sigma,\quad \delta\Omega= \dot\Sigma+ [\Omega,\Sigma], 
\la{redva} \ee
where $\Sigma(t)\in\mfrak{so}(3)$ specifies a variation vector field 
on $\mbf{SO}(3)$ that vanishes at the end-points, i.e., $\Sigma(0)=\Sigma(T)
=0$. Extremizing the action along this vector field, we get 
\beas \delta S= \int_0^T \left\{\frac12\langle\delta\Omega,\Omega\Lambda\rangle+ 
\frac12\langle\Omega,\delta\Omega\Lambda\rangle-\langle V_C,\Sigma\rangle\right\}
\di t \\ =\int_0^T \left\{\frac12\langle J(\Omega),\delta\Omega\rangle-
\langle \partial_C V,C\Sigma\rangle\right\}=0, \eeas
where $\partial_C V=\frac{\partial V}{\partial C}$ 
and $J:\mfrak{so}(3)\rightarrow\mfrak{so}(3)$ is a positive definite operator 
on the Lie algbera $\mfrak{so}(3)$ that is defined by 
\be J(\Omega)=\Lambda\Omega+\Omega\Lambda. \la{Jop} \ee 
The second term arising from the potential in the first variation above, can 
be rendered as 
\[ \langle C\T\partial_C V,\Sigma\rangle= \frac12 \langle C\T\partial_C V-
(\partial_C V)\T C,\Sigma\rangle. \]
Using the reduced variations as given by (\ref{redva}) and the relations 
(see \cite{miln}) 
\[ \langle J(A_1),A_2\rangle= \langle A_1,J(A_2)\rangle, \ \langle [X,Y],Z
\rangle= \langle X,[Y,Z]\rangle, \] 
we get (integrating by parts)
\beas 2\delta S &=& \int_0^T \left\{\langle J(\Omega),\dot\Sigma+[\Omega,
\Sigma]\rangle\right. \\ & &\left. -\langle C\T\partial_C V -(\partial_C V)\T C,
\Sigma\rangle\right\}\di t \\ &=& \int_0^T \left\{\langle [J(\Omega),\Omega],
\Sigma\rangle- \langle C\T\partial_C V-(\partial_C V)\T C,\Sigma\rangle 
\right. \\ & & \left. - \langle J(\dot\Omega),\Sigma
\rangle\right\}\di t +\langle J(\Omega),\Sigma\rangle\big|_0^T =0. \eeas
Since we take arbitrary fixed end-point variations, the last term above 
vanishes, and the terms in the integral give us the dynamics 
\be J(\dot\Omega)= [J(\Omega),\Omega]-C\T\partial_C V +(\partial_C V)\T C. 
\la{attdyn} \ee 
This dynamics equation is also derived in \cite{lelemc} in a similar fashion. 

We now present the result that any given value of $J(\Omega)$ for a given symmetric  
positive definite $\Lambda$, uniquely determines the skew-symmetric matrix $\Omega$. 
\begin{lemma}
If $K$ is symmetric and positive definite and $X\in \mfrak{so}(n)$, the map $J_K: 
\mfrak{so}(n)\rightarrow \mfrak{so}(n)$ given by $J_K: X \mapsto KX+XK$ has kernel 
zero, and is hence an isomorphism. 
\la{JXlem} 
\end{lemma}
\noindent{\em Proof}: Since $X$ is skew, there exists a unitary matrix $L$, i.e., 
$L\bar{L}\T=\bar{L}\T L= I_N$, such that $LX\bar{L}\T= \imath \Sigma$, where $\Sigma$ 
is a real diagonal matrix ($X$ is unitarily diagonalizable). Thus if $KX+XK=0$, then 
\[ LK\bar{L}\T\Sigma+ \Sigma LK\bar{L}\T=0, \]
where $LK\bar{L}\T$ is a positive definite Hermitian matrix. If $e$ is an eigenvector 
of $\widehat{K}= LK\bar{L}\T$, then 
\[ \widehat{K} e= \lambda e,\quad \lambda>0, \]
and 
\[ \widehat{K}\Sigma+ \Sigma\widehat{K}=0 \Rightarrow \widehat{K}(\Sigma e)+\lambda 
(\Sigma e)= 0. \]
Hence $\Sigma e$ is also an eigenvector of $\widehat{K}$ with eigenvalue $-\lambda 
<0$. But $\widehat{K}$ is positive definite and so all its eigenvalues are strictly 
positive. Thus, we have a contradiction, unless $\Sigma=0$ and hence $X=0$. \qed \\
This lemma and its proof are also given in \cite{san}. Note that we only need this 
result for $\mfrak{so}(3)$ here, although we have stated and proved it for all 
$\mfrak{so}(n)$. Thus, if the momentum $M=J(\Omega)$ is given, then one can obtain 
the unique angular velocity corresponding to it, $\Omega=J^{-1}(M)\in\mfrak{so}(3)$. 
From equation (\ref{attdyn}) and Lemma \ref{JXlem}, this uniquely determines 
$\dot\Omega$ given the values of $\Omega$ and $C$ at any instant. We use this 
dynamics equation, along with the attitude kinematics equation (\ref{atkin}), to 
propagate the attitude and angular velocity between discrete sets of measurements. 

\subsection{Attitude Estimation Filter without Angular Velocity Measurements}

We use the attitude determination algorithm given in Section \ref{sec:atder} 
to form an attitude estimation filter, by augmenting angular velocity data. The  
algorithm presented here works when there are body vector measurements at 
discrete time instants, but no available angular velocity measurements. However, 
we assume that we know the initial angular velocity. We use equations (\ref{atkin}) 
and (\ref{attdyn}) to integrate the attitude and angular velocity in time between 
the sets of body vector measurements; this corresponds to the propagation phase 
of the filter. The attitude and angular velocity are then updated based on the 
body vector measurements, in a manner similar to that used in a Kalman filter. 

We obtain the filter as an optimal filter from a suitable cost function that minimizes 
errors between the estimated and the measured attitude, as well as the difference 
between the propagated and updated estimates. The update of the attitude estimate 
$\wh{C_k}^+$ at measurement instant $t_k$ is obtained by minimizing the following 
cost function with respect to $\wh{C_k}^+$: 
\bea \cJ_a = \frac12 \sum_{k=0}^N \big\{\lan E_k-\wh{C_k}^+\wt{B_k}, (E_k-
\wh{C_k}^+\wt{B_k})W_k\ran \nn \\ +\lan (\wh{C_k}^-)\T\wh{C_k}^+ -I,((\wh{C_k}^-)\T
\wh{C_k}^+ -I)\Delta\ran \big\}, \la{cosfunc} \eea 
where $\Delta$ is a symmetric positive definite matrix that can be chosen as a design 
parameter for the filter. 
In the case that we are considering now, there are no angular velocity measurements; 
therefore we update the angular velocity by minimizing the difference between the 
rates of change of attitude at a measurement instant:
\be \cJ_r= \frac12 \lan \wh{C_k}^+\wh{\Omega_k}^+ - \wh{C_k}^-\wh{\Omega_k}^-,
(\wh{C_k}^+\wh{\Omega_k}^+ - \wh{C_k}^-\wh{\Omega_k}^-)\Pi\ran, \la{cfomeg} \ee 
where $\Pi$ is a symmetric positive definite matrix that can be chosen as a 
design parameter for the filter. Here $\wh{C_k}^-$ and $\wh{\Omega_k}^-$ are 
obtained by integrating equations (\ref{atkin}) and (\ref{attdyn}) respectively, 
from time $t_{k-1}$ to time $t_k$, with initial conditions $\wh{C_{k-1}}^+$ and 
$\wh{\Omega_{k-1}}^+$ respectively. 

The first variation of $\wh{C_k}^+$ is given by 
\[ \delta\wh{C_k}^+ = \wh{C_k}^+ U_k^+,\;\ U_k^+ \in\mfrak{so}(3). \]
Setting the first variation of $\cJ_a$ in (\ref{cosfunc}) to zero, we get: 
\bea \delta\cJ_a = \sum_{k=0}^N \big\{\lan (E_k-\wh{C_k}^+\wt{B_k})W_k\wt{B_k}\T,
-\wh{C_k}^+ U_k^+\ran \nn \\ +\lan (\wh{C_k}^+)\T\wh{C_k}^-\Delta\big(
(\wh{C_k}^-)\T \wh{C_k}^+ -I\big), U_k^+\ran\big\} \nn \\
= \sumk\big\{-\lan (\wh{C_k}^+)\T\big(\wh{C_k}^-\Delta+E_k W_k\wt{B_k}\T\big),
U_k^+\ran\big\}= 0, \la{extcosf} \eea
taking into account the initial condition $\wh{C_0}^+ = \wh{C_0}^-$,
which is either assumed to be known from a given initial attitude, or obtained from 
an initial set of measurements $\wt{B_0}$ using the algorithm of Theorem \ref{atde}. 
From the expression (\ref{extcosf}), we get the result 
\be (\wh{C_k}^+)\T L_k \mbox{  is symmetric, where  } L_k=\wh{C_k}^-\Delta+E_k W_k
\wt{B_k}\T. \la{exCkp} \ee 
Now the result of Theorem \ref{atde} can be applied to obtain the update of the 
attitude estimate $\wh{C_k}^+$ in terms of the QR decomposition of $L_k$. Given 
this update of the attitude estimate, one can obtain an update of the angular 
velocity estimate by minimizing $\cJ_r$ in (\ref{cfomeg}) with respect to 
$\wh{\Omega_k}^+$. This gives us:
\beas \delta\cJ_r &=& \lan\big(\wh{C_k}^+\wh{\Omega_k}^+ -\wh{C_k}^-\wh{\Omega_k}^-
\big)\Pi, \wh{C_k}^+\delta\wh{\Omega_k}^+\ran =0 \\ 
&\Rightarrow & \big(\wh{\Omega_k}^+ -(\wh{C_k}^+)\T\wh{C_k}^-\wh{\Omega_k}^-\big)
\Pi \mbox{  is symmetric.} \eeas
The initial angular velocity $\wh{\Omega_0}^-$ is assumed to be 
known, and from (\ref{inicon}), $\wh{\Omega_0}^- =\wh{\Omega_0}^+$. The above 
analysis can be formalized into the following result, which is one of the main 
results of this paper.   
\begin{theorem}
The attitude estimation filter obtained from minimizing the cost functions $\cJ_a$ 
and $\cJ_r$ is given by the attitude and angular velocity updates: 
\bea \wh{C_k}^+= S_k L_k,\quad \wh{\Omega_k}^+\Pi+ \Pi\wh{\Omega_k}^+ = (\wh{C_k}^+)\T
\wh{C_k}^-\wh{\Omega_k}^-\Pi \nn \\ +\Pi\wh{\Omega_k}^-(\wh{C_k}^-)\T\wh{C_k}^+, 
\la{filup} \eea  
where 
\be Q_k R_k =L_k =\wh{C_k}^-\Delta+E_k W_k\wt{B_k}\T, \la{LkQR} \ee 
is the QR decomposition of $L_k$, and 
\be S_k = Q_k\sqrt{(R_kR_k\T)^{-1}}Q_k\T \la{Skeq} \ee
is symmetric. The initial conditions for the filter are 
\be \wh{C_0}^+ = \wh{C_0}^-, \quad \wh{\Omega_0}^- =\wh{\Omega_0}^+, \la{inicon} \ee 
where $\wh{C_0}^-$ is either given or obtained from an initial set of measurements 
$\wt{B_0}$, and $\wh{\Omega_0}^-$ is given. The propagation equations for the filter 
are given by: 
\bea \wh{C_{k+1}}^- &=& \int_{t_k}^{t_{k+1}}\big(C\Omega\big)\di t, \la{atprop} \\ 
\wh{\Omega_{k+1}}^- &=& J^{-1}\left(\int_{t_k}^{t_{k+1}}\big([J(\Omega),\Omega]-
V_C(C)\big)\di t\right), \la{velpro} \eea 
where $C(t_k)=\wh{C_k}^+$ and $\Omega(t_k)= \wh{\Omega_k}^+$. 
\la{filt}
\end{theorem} 
Note that according to Lemma \ref{JXlem}, equations (\ref{filup}) and (\ref{velpro}) 
uniquely determine $\wh{\Omega_k}^+$ and $\wh{\Omega_{k+1}}^-$ respectively. 
The result of Proposition \ref{unbias} also holds, i.e., the attitude estimate given 
by this result is unbiased. This can be shown in a similar manner to the proof of 
Proposition \ref{unbias}. If $\Pi=I$ (the identity matrix), then equation (\ref{filup}) 
for the update of the angular velocity estimate simplifies to 
\be \wh{\Omega_k}^+ = \frac12 \left((\wh{C_k}^+)\T\wh{C_k}^-\wh{\Omega_k}^- 
+\wh{\Omega_k}^-(\wh{C_k}^-)\T\wh{C_k}^+\right). \la{anvup} \ee 
This equation can be readily used for updating angular velocity estimates 
without angular velocity measurements in the filter implementation. For the 
propagation equations (\ref{atprop}) and (\ref{velpro}) may be implemented 
by numerical integration software, including variational integrators 
(see \cite{haluwa, lelemc}) that preserve the group structure of $\mbf{SO}(3)$. 

\subsection{Attitude Estimation Filter using Angular Velocity Measurements} 

The creation of an attitude estimation filter from the basic attitude 
determination algorithm in Section \ref{sec:atder} is made easier when 
angular velocity measurements are available. In this case, we assume that 
the sampling instants for attitude and angular velocity measurements are the 
same. The body vector (attitude) measurements are given by (\ref{bomeas}), 
while the angular velocity measurements are given by 
\be \wt{\Omega_k}= \Omega_k + P_k, \la{angmes} \ee
where $\Omega_k=\Omega(t_k)$ is the actual angular velocity and $P_k$ is a 
zero mean measurement error from a stochastic process with known statistics. 
Extensions can also be made to deal with the case when the sampling 
instants for attitude and angular velocity measurements are different. 

The optimal filter when attitude and angular measurements are available is 
obtained by minimizing the following cost function with respect to $\wh{C_k}^+$ 
and $\wh{\Omega_k}^+$: 
\bea \cJ_b = \cJ_a+ \frac12 \sum_{k=0}^N \lan \wh{\Omega_k}^+ -\wt{\Omega_k},
(\wh{\Omega_k}^+-\wt{\Omega_k})X_k\ran, \la{cosfun1} \eea 
where $\cJ_a$ is as defined in (\ref{cosfunc}), and $\Delta$ and $\Gamma$ are 
symmetric positive definite matrices that can be chosen as design parameters 
for the filter. The matrix $X_k$ is also a symmetric positive definite matrix, 
which can be assigned as the error covariance matrix for the angular velocity 
measurement error $P_k$ in (\ref{angmes}). 

We take reduced variations on $\mbf{SO}(3)$ with the first variations of 
$\wh{C_k}^+$ and $\wh{\Omega_k}^+$ given by: 
\[ \delta\wh{C_k}^+ = \wh{C_k}^+ U_k^+,\;\ \delta\wh{\Omega_k}^+ =\dot U_k^+
+[\wh{\Omega_k}^+,U_k^+], \]
where $U_k^+ \in\mfrak{so}(3)$. The necessary condition for optimality is 
given by equating the first variation of $\cJ_b$ with respect to $\wh{C_k}^+$ 
and $\wh{\Omega_k}^+$ to zero. This gives us: 
\be \delta\cJ_b= \sum_{k=0}^N \big\{\lan [G_k,\wh{\Omega_k}^+]- (\wh{C_k}^+)\T
L_k,U_k^+ \ran+ \lan G_k,\dot U_k^+ \ran\big\}=0, \la{fstvar} \ee
where 
\bea G_k&=& (\wh{\Omega_k}^+-\wt{\Omega_k})X_k +(\wh{\Omega_k}^+ -\wh{\Omega_k}^-)
\Gamma, \nn \\ L_k &=& E_k W_k\wt{B_k}\T+\wh{C_k}^- \Delta. \la{kok} \eea 
Since $U_k^+$ and $\dot U_k^+$ are independent of each other, the second term 
in equation (\ref{fstvar}) implies that $G_k$ is symmetric. This in turn implies 
that $[G_k,\wh{\Omega_k}^+]$ is also symmetric, and hence $\lan [G_k,
\wh{\Omega_k}^+],U_k^+ \ran=0$. Thus, the first term in equation (\ref{fstvar}) 
implies that $(\wh{C_k}^+)\T L_k$ is symmetric. Now we can apply Theorem 
\ref{atde} to obtain the update of the attitude estimate. The attitide and 
angular velocity updates are given in the following result. 
\begin{theorem}
The attitude estimation filter obtained from minimizing the cost function $\cJ_a$ 
in equation (\ref{cosfun1}) is given by the attitude and angular velocity updates: 
\be \wh{C_k}^+= S_k L_k,\;\ J_{X_k+\Gamma}(\wh{\Omega_k}^+) = J_{X_k}
(\wt{\Omega_k})+ J_\Gamma (\wh{\Omega_k}^-), \la{filup1} \ee  
where $J_K: \mfrak{so}(n)\rightarrow \mfrak{so}(n)$  for a symmetric positive 
definite matrix $K$ is as defined in Lemma \ref{JXlem}, 
\be Q_k R_k =L_k =\wh{C_k}^-\Delta+E_k W_k\wt{B_k}\T, \la{LkQR1} \ee 
is the QR decomposition of $L_k$, and 
\be S_k = Q_k\sqrt{(R_kR_k\T)^{-1}}Q_k\T \la{Skeq1} \ee
is symmetric. The initial conditions for the filter are 
\be \wh{C_0}^+ = \wh{C_0}^-, \quad \wh{\Omega_0}^- =\wh{\Omega_0}^+, \la{incon1} \ee 
where $\wh{C_0}^-$ and $\wh{\Omega_0}^-$ are either given or obtained from initial 
measurements $\wt{B_0}$, and $\wt{\Omega_0}$. 
\la{filt1}
\end{theorem} 

The propagation equations for the filter are equations (\ref{atprop})-(\ref{velpro}) 
where $C(t_k)=\wh{C_k}^+$ and $\Omega(t_k)= \wh{\Omega_k}^+$. Thus, the propagation 
phases for the filters given by theorems \ref{filt} and \ref{filt1} are identical.  
Note that, by Lemma \ref{JXlem}, equation (\ref{filup1}) determines $\wh{\Omega_k}^+$ 
uniquely since $X_k+\Gamma$ is positive definite. Also, the result of Proposition 
\ref{unbias} holds for the attitude estimate given by this filter. The angular 
velocity estimate is also unbiased, since in the absence of angular velocity 
measurement errors, $\wt{\Omega_k}=\wh{\Omega_k}^-= \Omega_k$, and equation 
(\ref{filup1}) gives $\wh{\Omega_k}^+=\Omega_k$. The filters developed in this 
section can also be extended to estimate a constant bias in measurements, if sensor 
bias is present. 


\section{Conclusions}\la{sec:conc} 

This paper presents an attitude determination algorithm and attitude estimation 
filters that can be used for attitude estimation of robots, spacecraft, and other 
vehicles. The attitude determination algorithm is obtained from an optimization 
process with the cost function equal to a weighted attitude estimation error on 
the group of rigid-body orientations. This algorithm is global, and does not use 
any local coordinate representation (like Euler angles or quaternions) for the 
group of orientations. The optimization is carried out with variations on the smooth 
manifold (Lie group) of rigid body orientations, and the estimate obtained is 
shown to (globally) minimize the attitude estimation error. It is also shown to 
provide an unbiased estimate of the attitude, i.e., in the absence of measurement 
errors, the estimate of the attitude obtained is the actual attitude. A numerical 
simulation of this attitude determination algorithm, with a set of seven simulated 
body unit vector measurements with noise for seven given inertial unit vectors, is 
carried out. The order of error in the attitude estimate obtained using this 
algorithm is found to be no more than the order of the error in measurements. 

The attitude estimation filters are of the {\em continuous-discrete} type, which 
work with a continuous deterministic dynamics model supplemented by discrete sets 
of noisy measurements. These filters are obtained for two cases: when there are 
body vector measurements but no angular velocity measurements, and when there are 
attitude and angular velocity measurements at the same measurement instants. It 
is assumed that the attitude dynamics is deterministic, and accurately known. The 
(nonlinear) attitude kinematics and dynamics equations are used for propagation 
of the attitude and angular velocity between successive sets of measurements. The 
attitude and angular velocity estimates are obtained from an optimization process 
that minimizes the weighted sum of errors between the estimates and measurments, 
and between the estimates and propagated values for these quantities at each 
measurment instant. These filters are shown to be unbiased, i.e., in the absence 
of measurement errors, the estimates they give are equal to the actual attitude 
and angular velocity. 

This work is a preliminary exploration into attitude estimation techniques without 
using any local coordinate representation for the attitude. The results obtained 
thus far are encouraging, and they show that one can obtain unbiased filters that 
minimize the errors in attitude and angular velocity estimation without using 
local coordinate representations and without using an extended Kalman filter. While 
local coordinate representations of attitude have problems associated with 
singularities in attitude kinematics or additional constraints, the extended 
Kalman filter has problems associated with convergence of estimates for large 
initialization errors. These drawbacks are not present in the filters developed here,  
since they do not use local coordinates, and since they give optimal nonlinear 
filters that minimize the attitude and angular velocity estimation errors at each 
measurement instant. Thus, the attitude estimation filter algorithms developed here 
fill a gap in the existing research in this direction, besides improving upon the 
filters currently in use for attitude estimation of mechanical systems. 

Future work would include extension of the attitude and angular velocity 
estimation filters developed here to the case when the dynamics has modeling 
errors or noise. Numerical and/or experimental studies in implementation of 
these filters could also be explored. Numerical simulation results for the filters 
developed here, and numerical comparisons with estimation algorithms using local 
coordinates and extended Kalman filters for the deterministic dynamics case, have 
not been obtained yet. Such results are very likely to be reported in the near 
future.

\end{document}